\documentclass[a4paper]{amsart}

\usepackage{graphics,amssymb,enumerate}
\usepackage{hyperref}
\usepackage[dvips]{graphicx}
\usepackage[all]{xy}

\def\bd{\begin{displaymath}}
\def\ed{\end{displaymath}}
\newtheorem{theorem}{Theorem}
\newtheorem{prop}{Proposition}

\newtheorem{lemma}{Lemma}

\newtheorem{rem}{Remark}

\theoremstyle{definition}
\newtheorem{definition}{Definition}
\newtheorem{ex}{Example}

\newcommand{\mc}{\mathcal}
\newcommand{\p}{\partial}

\newcommand{\Rr}{\mathbb R}

\newcommand{\Zz}{\mathbb Z}

\renewcommand{\d}{\mathrm d}          
\renewcommand{\div}{\text{\rm div}}          

\newcommand{\X}{\ensuremath{\mathfrak{X}}}

\newcommand{\NN}{\ensuremath{\mathcal{N}}}
\newcommand{\Lie}{\boldsymbol{\pounds}}    
\newcommand{\tr}{\text{\rm tr}\,}          
\renewcommand{\mod}{\text{\rm mod}}          

\newcommand{\al}{\alpha}
\newcommand{\be}{\beta}

\newcommand{\la}{\lambda}

\Large
\begin{document}

\title{Integrable hierarchies and the modular class}

\author{Pantelis A.~Damianou}
\address{Department of Mathematics and Statistics\\
University of Cyprus\\
P.O.~Box 20537, 1678 Nicosia\\Cyprus}
\email{damianou@ucy.ac.cy}

\author{Rui Loja Fernandes}
\address{Departamento de Matem\'{a}tica\\
Instituto Superior T\'{e}cnico\\1049-001 Lisboa\\Portugal}
\email{rfern@math.ist.utl.pt}
\thanks{Supported in part by FCT/POCTI/FEDER and by grant
  POCTI/MAT/57888/2004.}

\begin{abstract}
We observe that the modular class of a Poisson-Nijhenhuis manifold has
a canonical representative and that, under a cohomological
assumption, this vector field is bi-hamiltonian. In many examples the
associated hierarchy of flows reproduces classical integrable hierarchies.
\end{abstract}

\date{July 2006}
\maketitle

\Large
\section{Introduction}
\label{intro}

It is well-known that the Poisson-Nijenhuis manifolds, introduced
by Kosmann-Schwarzbach and Magri in \cite{KosmannMagri}, form the
appropriate setting for studying  many classical integrable
hierarchies. In order to define the hierarchy, one usually
specifies \emph{in addition} to the Poisson-Nijenhuis manifold a
bi-hamiltonian vector field. In this paper we will show that to
every Poisson-Nijenhuis manifold one can associate a canonical
vector field (no extra choices are involved!) which under an
appropriate assumption defines an integrable
hierarchy of flows. Moreover, this vector field is a very natural
geometric entity, leading to a cohomological intrepertation of
this condition. For many classical examples we recover well-known
integrable hierarchies.

In order to explain in more detail our results, let us recall that
a Poisson manifold  $(M,\pi)$ usually does not carry a Liouville
form, i.e., a volume form which is invariant under the flows of
all hamiltonian vector fields(\footnote{We will assume that our
  manifolds are orientable. This is enough to cover all applications and
  simplifies the presentation. However, our results do extend to the
  non-orientable case.}).
The obstruction to the existence of an invariant volume form, as
was explained by J.-L.~Koszul \cite{koszul} and A.~Weinstein
\cite{weinstein2}, lies in the first Poisson cohomology group
$H^1_\pi(M)$ (the Poisson vector fields modulo hamiltonian vector
fields). More precisely, given a volume form $\mu$, we can
associate to it a Poisson vector field $X^{\pi}_\mu$, called the
\emph{modular vector field}. Though this vector field depends on
the choice of $\mu$, the Poisson cohomology class
$[X^{\pi}_\mu]\in H^1_\pi(M)$ does not, and this \emph{modular
class} is zero iff there exists some invariant measure on $M$. The
modular vector field was used by Dufour and Haraki in
\cite{dufour} to classify quadratic Poisson brackets in
$\mathbb{R}^3$. It was also useful in the classification of
Poisson structures in low dimensions, e.g., \cite{xu,marmo}.

Assume now that $(M,\pi_0,\mc{N})$ is a Poisson-Nijenhuis manifold
(\cite{KosmannMagri}). It is well
known that we can associate to it a hierarchy of Poisson structures:
\[ \pi_1:=\mc{N}\pi_0,\ \pi_2:=\mc{N}\pi_1=\mc{N}^2\pi_0,\dots\]
It is easy to check that the Nijenhuis tensor $\mc{N}$ maps
hamiltonian (respectively, Poisson) vector fields of $\pi_0$ to
hamiltonian (respect., Poisson) vector fields of $\pi_1$, and more
generally those of $\pi_i$ to those of $\pi_{i+1}$. However,
in general, for any choice of $\mu$, it \emph{does not} map the
modular vector field $X^0_\mu$ of $\pi_0$ to the modular vector
field $X^1_\mu$ of $\pi_1$. As we will show below, the difference:
\[ X_{\mc{N}}:=X^1_\mu-\mc{N}X^0_\mu,\]
is a Poisson vector field for $\pi_1$, which is \emph{independent}
of the choice of volume form $\mu$. Notice that this vector field
is zero if there exists a volume form $\mu$ which is invariant
simultaneously under the flows of the hamiltonian vector fields
for $\pi_1$ and $\pi_0$. Hence, we may think of $X_{\mc{N}}$ as a
modular vector field of our Poisson-Nijenhuis manifold. Moreover,
using the concept of \emph{relative} modular class, introduced
recently in \cite{GMM,KosmannWeinstein}, we can show that the
Poisson cohomology class of $X_{\mc{N}}$ is the relative modular
class of the transpose $\mc{N}^*$, when viewed as a morphism of
Lie algebroids. Furthermore, we will show the following result:

\begin{theorem}
Let $(M,\pi_0,\mc{N})$ be a Poisson-Nijenhuis manifold. Then the
modular vector field $X_{\mc{N}}$ is hamiltonian relative to
$\pi_0$ with hamiltonian equal to minus one half the trace of
$\mc{N}$:
\[ X_{\mc{N}}=X^0_{-\frac{1}{2}\tr \mc{N}}.\]
\end{theorem}

Therefore, the vector field $X_{\mc{N}}$ is hamiltonian relative
to $\pi_0$ and Poisson relative to $\pi_1$. So $X_{\mc{N}}$ is
very close to defining a bi-hamiltonian system, and hence  a
hierarchy of flows. Of course, the obstruction is the Poisson
cohomology class $[X_{\mc{N}}]\in H^1_{\pi_1}(M)$, i.e., the modular class of
the Poisson-Nijenhuis manifold. This
class is zero, for example, if there are measures $\mu$ and $\eta$
invariant under both the hamiltonian flows of $\pi_0$ and $\pi_1$.
Note that, in general, $X_{\mc{N}}$ itself will still be non-zero,
in which case the two invariant measures are non-proportional. A
typical situation that fits many examples is the following:

\begin{theorem} \label{trace}
Let $(M,\pi_0,\mc{N})$ be a Poisson-Nijenhuis manifold and assume
that $\mc{N}$ is non-degenerate. Then the
modular vector field $X_{\mc{N}}$ is bi-hamiltonian and hence
determines a hierarchy of flows which are given by:
\[ X_{i+j}=\pi_i^\sharp\d h_j=\pi_j^\sharp\d h_i\quad (i,j \in {\bf Z} )\]
where
\[ h_0=-\frac{1}{2}\,\log(\det\mc{N}),\quad
h_i=-\frac{1}{2i}\,\tr\mc{N}^i \quad (i \not=0 ).\]
\end{theorem}

We will see below that most of the known hierarchies of integrable
systems can be obtained in this manner, therefore providing a new
approach to the integrability of those systems. Moreover, in some
cases (e.g., the Toda systems) it gives rise to previously unknown
bi-hamiltonian formulations. Note that the fact that the traces of the
powers of $\mc{N}$ give rise to a hierarchy of flows was noticed early
in the history of integrable systems (see, e.g., \cite{Kosmann,Magri}).

The paper is organized as follows. In Section 2, we recall a few
basic facts concerning modular vector fields and modular classes,
and we show that the modular class of the Lie algebroid associated
with a Nijenhuis tensor is represented by $\d(\tr\mc{N})$. This
basic fact, which does not seem to have been noticed before, sets
up the stage for section 3, where we consider the modular vector
field of a Poisson-Nijenhuis manifold. In section 4, we introduce
integrable hierarchies related to the modular class and we prove
Theorem \ref{trace} above. In Section 5, we show how one can
recover many of the known classical integrable hierarchies using
our results.

\textbf{Acknowledgments.} We would like to thank several institutions
for their hospitality while work on this project was being done:
Instituto Superior T\'ecnico and Universit\'e de Poitiers (Pantelis Damianou);
University of Cyprus, University of Milano-Bicoca and ESI Vienna
(Rui L.~Fernandes). We would like to thank Yvette
Kosmann-Schwarzbach for many comments on a first version of
this paper, which helped improving it greatly, Franco Magri
who pointed out to us that the assumptiom (made on the same
first version of the paper) of invertibility of $\pi_0$ is actually
superfulous, and Raquel Caseiro for usefull discussions.

\section{Modular classes}

In this section we present several results concerning modular
classes that will be needed later. This will also help in
establishing our notation. Our main result here is Proposition
\ref{prop:mod:class:N}, where we compute the modular class of the
Lie algebroid associated with a (1,1)-tensor $\mc{N}$ with
vanishing Nijenhuis torsion.

\subsection{Modular class of a Poisson manifold}
If $(M,\{\cdot,\cdot\})$ is a Poisson manifold, we will denote by
$\pi\in\X^2(M)$ the associated Poisson tensor which is given by
\[ \pi(\d f,\d  g):=\{f,g\}, \quad (f,g \in  C^\infty(M))\]
and by $\pi^\sharp:T^\ast M \rightarrow TM$ the vector bundle map defined by
\[
\pi^\sharp (\d h)=X_h:=\{h,\cdot\},
\]
where $X_h$ is the hamiltonian vector field determined by $h\in
C^\infty(M)$. Recall also that the Poisson cohomology of
$(M,\pi)$, introduced by Lichnerowicz \cite{lich}, is the
cohomology of the complex of multivector fields
$(\X^\bullet(M),\d_\pi)$, where the coboundary operator is defined
by taking the Schouten bracket with the Poisson tensor:
\[ \d_{\pi}A\equiv [\pi, A].\]
This cohomology is denoted  by $H^\bullet_\pi(M)$. We will be
mainly interested in the first  Poisson cohomology space
$H^1_\pi(M)$, which is just the space of Poisson vector fields
modulo the hamiltonian vector fields. Note that our conventions
are such that the hamiltonian vector field associated with the
function $h$ is given by:
\begin{equation}
\label{eq:hamiltons} X_h=-[\pi,h]=-\d_\pi h.
\end{equation}
In this paper we follow the same sign conventions as in the book by
Dufour and Zung \cite{DufourZung}, and which differ from other sign
conventions such as the one in Vaisman's monograph \cite{vaisman}
(\footnote{In particular the Schouten bracket on
 multivector fields satisfies the following super-commutation,
 super-derivation and super-Jacobi identities:
 \begin{align*}
     &[A,B]=-(-1)^{(a-1)(b-1)} [B,A]\\
     &[A,B\wedge C]=[A,B]\wedge C+(-1)^{(a-1)b}B\wedge[A,C]\\
     &(-1)^{(a-1)(c-1)} [A,[B,C]]+ (-1)^{(b-1)(a-1)} [B,[C,A]] +(-1)^{(c-1)(b-1)} [C,[A,B]]=0
 \end{align*}
 where $A\in\X^a(M)$, $B\in\X^b(M)$ and $C\in\X^c(M)$.}).

Let us assume that $M$ is oriented and fix an arbitrary volume form
$\mu\in\Omega^{\text{top}}(M)$. The divergence of a vector field
$X\in\X(M)$ relative to $\mu$ is the unique function $\div_\mu(X)$
that satisfies:
\[ \Lie_X\mu=\div_\mu(X)\mu.\]
When $(M,\pi)$ is a Poisson manifold, a volume form $\mu$ defines
the modular vector field:
\[ X_\mu(f):=\div_\mu(X_f).\]
Note that this vector field depends on the choice of $\mu$.

More generally, a choice of volume form $\mu$ induces, by contraction, an
isomorphism $\Phi_\mu:\X^k(M)\to\Omega^{m-k}(M)$, where $m=\dim M$,
and we define, following Koszul \cite{koszul}, the following operator
that generalizes the divergence operator above:
$D_\mu:\X^k(M)\to\X^{k-1}(M)$ defined by:
\[ D_\mu=\Phi_\mu^{-1} \circ \d \circ \Phi_\mu \ , \]
where $\d$ is the exterior derivative. It is obvious that $D_\mu^2=0$,
so $D_\mu$ is a homological operator. Now we have:

\begin{prop}
\label{prop:mod:Poisson}
For a Poisson manifold $(M,\pi)$ with a  volume  form $\mu$ the
modular vector field is given by:
\begin{equation}
  \label{eq:modular:vf}
  X_\mu=D_\mu(\pi).
\end{equation}
If $(x^1,\dots,x^m)$ are local coordinates, such that $\mu=\d
x^1\wedge\cdots\wedge\d x^m$ and $\pi=\sum_{i<j}\pi^{ij}
\frac{\partial}{\partial x^i}\wedge\frac{\partial}{\partial x^j}$
then:
\[
X_\mu=\sum_{i=1}^{m}\left( \sum_{j=1}^m \frac{\partial
\pi^{ij}}{\partial x_j} \right)\frac{\partial}{\partial x_i}.
\]
\end{prop}

The proof of this proposition is standard and we refer, for example, to
\cite[Chapter 2.6]{DufourZung} for details.

Some authors take expression (\ref{eq:modular:vf}) as the definition
of the modular vector field. Recalling that the Koszul operator
satisfies the basic identity:
\begin{equation}
  \label{eq:Koszul:Schouten}
  D_\mu([A,B])=[A,D_\mu(B)]+(-1)^{b-1}[D_\mu(A),B]
\end{equation}
we see immediately from $[\pi,\pi]=0$ that
\[ \d_\pi X_\mu=[\pi,D_\mu(\pi)]=0,\]
so the modular vector field is a Poisson vector field. Also, if we
are given another volume form $\mu'$, so that $\mu'=g\mu$ for some
non-vanishing function $g$, we find from the definition of the
Koszul operator:
\[
  D_{g\mu}A=D_\mu A+[A,\ln|g|].
\]
In particular, when $A=\pi$ this shows that under a change of
volume form the modular vector field changes by an addition of a
hamiltonian vector field:
\begin{equation}
\label{eq:change:form} X_{g\mu}=X_\mu-X_{\ln|g|}.
\end{equation}
Therefore, the class $\mod(\pi)\equiv[X_\mu]\in H^1_\pi(M)$ is well-defined.

\subsection{Modular class  of a Lie algebroid}

We will need also the modular class of a Lie algebroid, which was
introduced in \cite{EvLuWe}.

Let $p:A\to M$ be a Lie algebroid over $M$, with anchor $\rho:A\to TM$
and Lie bracket
$[\cdot,\cdot]:\Gamma(A)\times\Gamma(A)\to\Gamma(A)$. Lie algebroids
are some kind of generalized tangent bundles, so many of the
constructions from the usual tensor calculus can be extended to Lie
algebroids, and we recall a few of them. First, the algebroid
cohomology of $A$ is the cohomology of the complex
$(\Omega^k(A),\d_A)$, where $\Omega^k(A)\equiv\Gamma(\wedge^k A^*)$
and $\d_A:\Omega^k(A)\to \Omega^{k+1}(A)$ is the de Rham type
differential:
\begin{multline}
\label{eq:A:differential} \d_A \omega(\al_0,\dots,\al_k)=
\sum_{i=0}^k (-1)^{i+1}
\rho(\al_i)(\omega(\al_0,\dots,\check{\al}_i,\dots,\al_k))\\
\sum_{0\le i<j\le k} (-1)^{i+j}
\omega([\al_i,\al_j],\al_0,\dots,\check{\al}_i,\dots,\check{\al}_j,\dots,\al_k)
\ .
\end{multline}
The Lie algebroid cohomology is denoted by $H^\bullet(A)$. Given a
section $\al\in \Gamma(A)$ (a ``vector field''), there is a Lie
$A$-derivative operator $\Lie_\al$ and a contraction operator $i_\al$
defined as in the usual case of $TM$, but using the $A$-Lie bracket. It
follows that we also have Cartan's magic formula (for details see,
e.g., \cite{EvLuWe}):
\[ \Lie_\al=i_\al\d_A+d_A i_\al.\]

Now to define the modular class of $A$ we proceed as follows. We
assume that the line bundles $\wedge^{\text{top}}A$ and
$\wedge^{\text{top}}T^*M$ are trivial and we choose global
sections $\eta$ and $\mu$   (\footnote{Again, this
  orientability assumption is made only to simplify the presentation,
  and is not essential for what follows.}).
Then $\eta\otimes\mu$ is a section of
$\wedge^{\text{top}}A\otimes\wedge^{\text{top}}T^*M$, and we define
$\xi_A\in C^1(A)$ to be the unique element such that:
\[
(\Lie_\al\eta)\otimes\mu+\eta\otimes(\Lie_{\rho(\al)}\mu)
=\xi_A(\al)\eta\otimes\mu, \quad \forall \al\in\Gamma(A).
\]
One checks that $\xi_A$ is indeed an $A$-cocycle, and that its
cohomology class is independent of the choice of $\eta$ and
$\mu$. Hence, there is a well-defined modular class of $A$ denoted
$\mod(A)\equiv [\xi_A]\in H^1(A)$.

\begin{ex}
\label{ex:mod:Poisson}
For a Poisson manifold $(M,\pi)$ we have a natural Lie algebroid
structure on its cotangent bundle $T^*M$. For the anchor we have
$\rho=\pi^\sharp$ and for the Lie bracket on sections of $A=T^*M$, i.e.,
on one forms, we have:
\[ [\al,\be]=\Lie_{\pi^\sharp\al}\be-\Lie_{\pi^\sharp\be}\al-\d\pi(\al,\be).\]
Note however that the two definitions above of the modular class differ
by a multiplicative factor:
\[ \mod(T^*M)=2~\mod(\pi).\]
See also \cite{fernandes2} for a deeper explanation of the factor
2. This factor will appear frequently in our formulas.
\end{ex}

\begin{ex}
\label{ex:mod:Nijenhuis}
Let $M$ be a manifold and $\mc{N}:TM\to TM$ a Nijenhuis tensor,
i.e., a (1,1)-tensor whose Nijenhuis torsion
\[ T_{\mc{N}}(X,Y):=\mc{N}[\mc{N}X,Y]+\mc{N}[X,\mc{N}Y]-\mc{N}^2([X,Y])-[\mc{N}X,\mc{N}Y],\]
vanishes. This is equivalent to requiring that the triple
$(TM,[~,~]_\mc{N},\rho_\mc{N})$ is  a Lie algebroid, where the
anchor is given by
\[ \rho_\mc{N}(X):=\mc{N}X,\]
and the Lie bracket is defined by:
\[ [X,Y]_\mc{N}:=[\mc{N}X,Y]+[X,\mc{N}Y]-\mc{N}([X,Y]).\]
Let us compute the modular class of this Lie algebroid.

\begin{prop}
\label{prop:mod:class:N}
The modular class of $(TM,[~,~]_\mc{N},\rho_\mc{N})$ is the cohomology
class represented by the 1-form $\d(\tr\mc{N})$.
\end{prop}

Note that this class may not be trivial: we must consider it as a
cohomology class in the Lie algebroid cohomology of
$(TM,[~,~]_\mc{N},\rho_\mc{N})$. This cohomology is computed by
the complex of differential forms but with a modified differential
$\d_\mc{N}$ that satisfies:
\[ \d_\mc{N}\mc{N}^*=\d\mc{N}^* \]
(here $\mc{N}^*:T^*M\to T^*M$ denotes the transpose of $\mc{N}$).

\begin{proof}[Proof of Proposition \ref{prop:mod:class:N}]
We pick a volume form $\mu\in\Omega^{\text{top}}(M)$, and we let
$\eta\in\Gamma(\wedge^{\text{top}}A)=\X^{\text{top}}(M)$ be the dual
multivector field: $\langle \mu,\eta\rangle=1$. Around any point, we
can choose local coordinates $(x^1,\dots,x^m)$ such that:
\[ \mu=\d x^1\wedge\cdots\wedge\d x^m,\quad
\eta=\frac{\partial}{\partial x^1}\wedge\cdots\wedge
\frac{\partial}{\partial x^m} .\]
In these coordinates, we write
\[ \mc{N}=\sum_{i,j=1}^m N^i_j \frac{\partial}{\partial x^i}\otimes\d x^j,\]
and for $X=\frac{\partial}{\partial x^k}$ we compute:
\begin{align*}
\Lie^\mc{N}_X \frac{\partial}{\partial x^i}&=
\sum_{j=1}^m
\left(\frac{\partial N^j_i}{\partial x^k}-\frac{\partial
    N^j_k}{\partial x^i}\right)\frac{\partial}{\partial x^j}\\
\Lie_{\mc{N}X}\d x^i&=\sum_{j=1}^m \frac{\partial N^i_k}{\partial
x^j}\d x^j \ ,
\end{align*}
where the first Lie derivative is in the Lie algebroid sense, while
the second is the usual Lie derivative. From these expressions it
follows that:
\begin{align*}
\Lie^\mc{N}_X\eta&=
\sum_{i=1}^m \frac{\partial}{\partial x^1}\wedge\cdots\wedge
\Lie^\mc{N}_X \frac{\partial}{\partial
  x^i}\wedge\cdots\wedge\frac{\partial}{\partial x^m}\\
&=\sum_{i=1}^m
\left(\frac{\partial N^i_i}{\partial x^k}-\frac{\partial
    N^i_k}{\partial x^i}\right)\eta
\end{align*}
and, similarly, that:
\begin{align*}
\Lie_{\mc{N}X}\mu
&=\sum_{i=1}^m \d x^1\wedge\cdots\wedge \Lie_{\mc{N}X}\d x^i \wedge\cdots\wedge\d x^m\\
&=\sum_{i=1}^m \frac{\partial N^i_k}{\partial x^i}\mu \ .
\end{align*}
Therefore, we conclude that for $X=\frac{\partial}{\partial x^k}$:
\[
\Lie^\mc{N}_X\eta\otimes\mu+\eta\otimes\Lie_{\mc{N}X}\mu
=\sum_{i=1}^m \frac{\partial N^i_i}{\partial x^k}\eta\otimes\mu
=\langle \d(\tr \NN),X\rangle \eta\otimes\mu.
\]
By linearity, this formula holds for every vector field $X$, on
any coordinate neighborhood. Hence, it must hold on all of $M$. We
conclude that $\d(\tr \NN)$ represents the modular class of
$(TM,[~,~]_\mc{N},\rho_\mc{N})$.
\end{proof}

Note that we have chosen the volume forms $\eta=\mu^{-1}$. For other
choices of $\eta$ and $\mu$ we would obtain different representatives
of the modular class. However, whenever we choose
$\eta=\mu^{-1}$ we always get the same representative, independent of
the choice of $\mu$. Also, the appearance of the trace should not be a
surprise in view of the interpretation of the modular class as a
secondary characteristic class (see \cite{fernandes2,fernandes3})
associated with the trace.
\end{ex}

\subsection{Relative modular class}
Let $\phi:A\to B$ be a morphism of Lie algebroids over the
identity. Then we have an induced chain map
$\phi^*:(\Omega^k(B),d_B)\to(\Omega^k(A),d_A)$ defined by:
\[ \phi^*P(\al_1,\dots,\al_k)=P(\phi(\al_1),\dots,\phi(\al_k)), \]
and, hence, also a morphism at the level of cohomology:
\[ \phi^*:H^k(B)\to H^k(A).\]

We can attach to this morphism a \emph{relative modular class}. Again, we
assume that $\wedge^{\text{top}}A$ and $\wedge^{\text{top}}B$ are
trivial line bundles, so we take global sections
$\eta\in\Gamma(\wedge^{\text{top}}A)$ and
$\nu\in\Gamma(\wedge^{\text{top}}B^*)$. Then we can define
$\xi_{A,B}^\phi\in C^1(A)$ to be the unique element such that:
\[
(\Lie^A_\al\eta)\otimes\mu+\eta\otimes(\Lie^B_{\phi(\al)}\mu)
=\xi_{A,B}^\phi(\al)\eta\otimes\mu, \quad \forall \al\in\Gamma(A).
\]
One can check that $\xi_{A,B}^\phi$ is in fact a cocycle, and that
its cohomology is independent of the choice of trivializing
sections $\eta$ and $\nu$. We conclude that we have a well defined
relative modular class:
\[\mod(A,B,\phi)\equiv \xi_{A,B}^\phi\in H^1(A).\]

Now we have the following basic fact (see \cite{KosmannWeinstein,GMM}):

\begin{prop}
\label{prop:mod:morph}
Let $\phi:A\to B$ be a morphism of Lie algebroids. Then:
\begin{equation}
\label{eq:rel:modular:class}
\mod(A,B,\phi)=\mod(A)-\phi^*\mod(B).
\end{equation}
Moreover, if $\psi:B\to C$ is another morphism, we have
\begin{equation}
\label{eq:compos:modular:class}
\mod(A,C,\psi\circ\phi)=\mod(A,B,\phi)+\phi^*\mod(B,C,\psi).
\end{equation}
\end{prop}

If we make any choice of sections
$\eta\in\Gamma(\wedge^{\text{top}}A)$,
$\nu\in\Gamma(\wedge^{\text{top}}B)$,
$\mu\in\Gamma(\wedge^{\text{top}}T^*M)$, and we choose
$\nu'\in\Gamma(\wedge^{\text{top}}B^*)$ to be dual to $\nu$,
(i.e., $\langle \nu,\nu'\rangle=1$), then
(\ref{eq:rel:modular:class}) already holds at the level of
cocycles, not just of cohomology classes: in the notation above,
we have the equality
\[ \xi_{A,B}^\phi=\xi_A-\phi^*\xi_B.\]
Similarly, (\ref{eq:compos:modular:class}) is also true at the level
of cocycles.

\begin{ex}
The tangent bundle $TM$ of any manifold is a Lie algebroid for the
usual Lie bracket of vector fields and the identity map as an
anchor. For this Lie algebroid, if we take a section $\nu\in
\Gamma(\wedge^{\text{top}}TM)$ and its dual section $\mu\in
\Gamma(\wedge^{\text{top}}TM^*)$, we see immediately that
$\xi_{TM}=0$, so its modular class vanishes. Now, given any Lie
algebroid $(A,[\cdot,\cdot])$, its anchor $\rho:A\to TM$ is a Lie
algebroid morphism. Hence, we conclude that
\[ \mod(A,TM,\rho)=\mod(A).\]
In particular, in the case of a Poisson manifold $(M,\pi)$ we find:
\[ \mod(T^*M,TM,\pi^\sharp)=\mod(T^*M)=2\,\mod(\pi).\]
Again, this equality is true already at the level of vector fields.
\end{ex}

\section{Modular vector fields and Poison-Nijenhuis manifolds}

We are now ready to look at Poisson-Nijenhuis manifolds and their
modular classes.

\subsection{Poisson-Nijhenhuis manifolds}
Let $(M,\pi_0,\mc{N})$ be a Poisson-Nijhenhuis manifold. Let us recall
what this means (\cite{KosmannMagri}):
\begin{enumerate}[(i)]
\item $\pi_0$ is a Poisson structure on $M$;
\item $\mc{N}:TM\to TM$ is a Nijenhuis tensor;
\item $\pi_0$ and $\mc{N}$ are compatible.
\end{enumerate}
The compatibility of $\pi_0$ and $\mc{N}$ means, first of all, that
\begin{equation}
    \label{eq:PN:compatibiliy:1}
    \mc{N}\pi_0^\sharp=\pi_0^\sharp\mc{N}^*,
\end{equation}
so that $\pi_1=\mc{N}\pi_0$ is a bivector field, and secondly
that the bracket on 1-forms $[~,~]_{\pi_1}$ naturaly associated
with $\pi_1$ (see Example \ref{ex:mod:Poisson}):
\[ [\al,\be]_{\pi_1}:=\Lie_{\pi_1^\sharp\al}\be-\Lie_{\pi_1^\sharp\be}\al-\d\pi_1(\al,\be)\]
and the bracket $[~,~]_{\pi_0}^{\mc{N}^*}$ obtained from $[~,~]_{\pi_0}$ by twisting by $\mc{N}^*$
(see Example \ref{ex:mod:Nijenhuis}):
\[ [\al,\be]_{\pi_0}^{\mc{N}^*}:=
[\mc{N}^*\al,\be]_{\pi_0}+[\al,\mc{N}^*\be]_{\pi_0}-\mc{N}^*([\al,\be]_{\pi_0})\]
actually coincide:
\begin{equation}
    \label{eq:PN:compatibiliy:2}
    [\al,\be]_{\pi_1}=[\al,\be]_{\pi_0}^{\mc{N}^*}.
\end{equation}

As a consequence of this definition, we have that $\pi_1$ must be a
Poisson tensor and the dual of the Nijenhuis tensor:
\[
\mc{N}^*:(T^*M,[~,~]_{\pi_1},\pi_1^\sharp)\to (T^*M,[~,~]_{\pi_0},\pi_0^\sharp)
\]
is a morphism of Lie algebroids.

As is well-known (\cite{KosmannMagri}), we have in fact a whole hierarchy
of Poisson structures:
\[ \pi_1:=\mc{N}\pi_0,\ \pi_2:=\mc{N}\pi_1=\mc{N}^2\pi_0,\dots\]
which are pairwise compatible:
\[ [\pi_i,\pi_j]=0,\quad\forall i,j=0,1,2,\dots\]
From this it follows that if we have a bi-hamiltonian vector field:
\[ X_1=\pi_1^\sharp\d h_0=\pi_0^\sharp\d h_1,\]
then we have a whole hierarchy of commuting flows $X_1,X_2,X_3,\dots$
where the higher order flows are given by:
\[ X_i=\pi_i^\sharp\d h_0=\pi_{i-1}^\sharp\d h_1.\]
Hence, one usually thinks of an integrable hierarchy as being
specified by a Poisson-Nijenhuis manifold \emph{and} a
bi-hamiltonian vector field. Here we would like to show that,
under a \emph{natural} assumption, there is a \emph{canonical hierarchy}
associated with a Poisson-Nijenhuis manifold, which does not
involve other choices such as a specification of a bi-hamiltonian
vector field. The source of this hierarchy is the modular class of
a Poisson-Nijenhuis manifold.

\subsection{Modular vector field of a Poisson-Nijhenhuis manifold}

Let $(M,\pi_0,\mc{N})$ be a Poisson-Nijhenhuis manifold. It is clear
from the definition that $\mc{N}$ maps the hamiltonian vector field
$X^0_f$ (relative to $\pi_0$) to the hamiltonian vector field $X^1_f$
(relative to $\pi_1$). Similarly, it is easy to see that $\mc{N}$ maps
Poisson vector fields of $\pi_0$ to Poisson vector fields of $\pi_1$.
More generally, $\mc{N}$ induces a map at the level of multivector
fields, denoted by the same letter
$\mc{N}:\X^\bullet(M)\to\X^\bullet(M)$, which is defined by:
\[ \mc{N}A(\al_1,\dots,\al_a)=A(\mc{N}^*\al_1,\dots,\mc{N}^*\al_a).\]
We have:

\begin{prop}
\label{prop:recursion:cohomology}
The map $\mc{N}:(\X^\bullet(M),\d_{\pi_0})\to (\X^\bullet(M),\d_{\pi_1})$
is a morphism of complexes:
\[ \mc{N}\d_{\pi_0}=\d_{\pi_1}\mc{N}.\]
\end{prop}

\begin{proof}
We need simply to observe that we have a Lie algebroid morphism:
\[
\mc{N}^*:(T^*M,[~,~]_{\pi_1},\pi_1^\sharp)\to (T^*M,[~,~]_{\pi_0},\pi_0^\sharp)
\]
so it induces a morphism between the complexes of forms of these Lie
algebroids, in the opposite direction. Of course, this map is just the map
$\mc{N}:(\X^\bullet(M),\d_{\pi_0})\to (\X^\bullet(M),\d_{\pi_1})$
introduced above.
\end{proof}

It follows that we have an induced map in cohomology
\[ \mc{N}:H^\bullet_{\pi_0}(M)\to H^\bullet_{\pi_1}(M).\]
Note, however, that in general $\mc{N}$ \emph{does not} map the
modular class of $\pi_0$ to the modular class of $\pi_1$. For a
choice of volume form $\mu\in\Omega^{\text{top}}(M)$, let us
denote denote by $X^1_\mu$ and by $X^0_\mu$ the modular vector
fields associated with $\pi_1$ and $\pi_0$ respectively.

\begin{lemma}
If $\mu$ and $\mu'$ are any two volume forms then:
\[ X^1_\mu-\mc{N}X^0_\mu=X^1_{\mu'}-\mc{N}X^0_{\mu'}.\]
Moreover, this vector field is Poisson relative to $\pi_1$.
\end{lemma}

\begin{proof}
Let $g\in C^\infty(M)$ be a non-vanishing function such that
$\mu'=g\mu$. By relation (\ref{eq:change:form}), we have
\begin{align*}
X^1_{\mu'}-\mc{N}X^0_{\mu'}
&=X^1_{g\mu}-\mc{N}X^0_{g\mu}\\
&=X^1_{\mu}-X^1_{\ln|g|}-\mc{N}(X^0_{\mu}-X^0_{\ln|g|})\\
&=X^1_\mu-\mc{N}X^0_\mu,
\end{align*}
where we used that $X^1_f=\mc{N}X^0_f$, for any function $f$.

The modular vector field $X^1_\mu$ is a Poisson vector field relative
to $\pi_1$. On the other hand, $\mc{N}$ maps the vector field
$X^0_\mu$, which is Poisson relative to $\pi_0$, to a Poisson vector
field relative to $\pi_1$. Hence, the sum $X^1_\mu-\mc{N}X^0_\mu$ is a
Poisson vector field relative to $\pi_1$.
\end{proof}

Let us set $X_\mc{N}:=X^1_\mu-\mc{N}X^0_\mu$ (which, by the lemma, is
independent of the choice of $\mu$). Note that
$X_\mc{N}$ is a vector field intrinsically associated with the
Poisson-Nijenhuis manifold $(M,\pi_0,\mc{N})$.

\begin{definition}
The vector field $X_\mc{N}$ is called the \textbf{modular vector field of
the Poisson-Nijenhuis} manifold $(M,\pi_0,\mc{N})$.
\end{definition}

The modular vector field $X_\mc{N}$ of $(M,\pi_0,\mc{N})$ will
play a fundamental role in the sequel. Our next proposition gives
further justification for this name and explains the possible
failure in $X_\mc{N}$ being a hamiltonian vector field:

\begin{prop}
\label{prop:rel:mod:PN}
Let $(M,\pi_0,\mc{N})$ be a Poisson-Nijenhuis manifold. The Poisson
cohomology class $[X_\mc{N}]\in H^1_{\pi_1}(M)$ equals half the
relative modular class of the Lie algebroid morphism:
\[
\mc{N}^*:(T^*M,[~,~]_{\pi_1},\pi_1^\sharp)\to (T^*M,[~,~]_{\pi_0},\pi_0^\sharp)
\]
\end{prop}

\begin{proof}
By Proposition \ref{prop:mod:morph} and Example \ref{ex:mod:Poisson}
we find:
\begin{align*}
\mod(T^*M_{\pi_1},T^*M_{\pi_0},\mc{N}^*)
&=\mod(T^*M_{\pi_1})-(\mc{N}^*)^*\mod(T^*M_{\pi_0})\\
&=2\,\mod(\pi_1)-2\,\mc{N}\mod(\pi_0)\\
&=2\,[X^1_\mu]-2\,\mc{N}[X^0_\mu]\\
&=2\,[X^1_\mu-\mc{N}X^0_\mu]=2\,[X_\mc{N}],
\end{align*}
for any volume form $\mu$.
\end{proof}

We emphasize that $X_\mc{N}$ is a canonical representative of the
relative modular class of $\mc{N}^*$, which does not depend on any
choice of measure.

\subsection{Hamiltonian character of the modular vector field}
As we saw above, the modular vector field
$X_\mc{N}$ of a Poisson-Nijenhuis manifold $(M,\pi_0,\mc{N})$ is a
Poisson vector field relative to $\pi_1$, which may fail to be
hamiltonian. Let us now look at its behavior relative to $\pi_0$.
We have:

\begin{theorem}
\label{thm:main:1} Let $(M,\pi_0,\mc{N})$ be a Poisson-Nijenhuis
manifold. Then the modular vector field $X_{\mc{N}}$ is hamiltonian
relative to $\pi_0$ with hamiltonian equal to minus one half the trace of
$\mc{N}$:
\begin{equation}
\label{eq:mod:0}
X_{\mc{N}}=X^0_{-\frac{1}{2}\tr \mc{N}}.
\end{equation}
\end{theorem}

Before we prove this theorem, let us  observe that this result is
intimately related to Proposition \ref{prop:mod:class:N}, where we
showed that the modular class of the Lie algebroid of a Nijenhuis
tensor $\mc{N}$ is represented by the 1-form $\d(\tr\mc{N})$. In
fact, observe that the compatibility condition of a Poisson-Nijenhuis
structure states that the two Lie algebroids
\begin{align*}
T^*M_{\pi_1}&:=(T^*M,[~,~]_{\mc{N}\pi_0},\rho=\mc{N}\pi_0^\sharp)\\
T^*M_{\pi_0}^{\mc{N}^*}&:=(T^*M,[~,~]_{\pi_0}^{\mc{N}^*},\rho=\pi_0^\sharp\mc{N}^*)
\end{align*}
actually coincide. Therefore they have the same modular classes, and
from the general Lie algebroid version of Proposition \ref{prop:mod:class:N}
we obtain:
\[ \mod(T^*M_{\pi_1})=\mod(T^*M_{\pi_0}^{\mc{N}^*})=[\d_{\pi_0}(\tr\mc{N})]+
\mc{N}^*\mod(T^*M_{\pi_0}).\]
Using Proposition \ref{prop:mod:morph}, this leads immediately to the statement:
\begin{align*}
2[X_{\mc{N}}]&=\mod(T^*M_{\pi_1},T^*M_{\pi_0},\mc{N}^*)\\
&=\mod(T^*M_{\pi_1})-\mc{N}^*\mod(T^*M_{\pi_0})\\
&=[\d_{\pi_0}(\tr\mc{N})]=-[(\pi_0^\sharp)\d(\tr\mc{N})].
\end{align*}
By working at the level of representatives of these cohomology
classes, one can give a proof of Theorem \ref{thm:main:1}. However, we
prefer to give a local coordinate proof which is a direct translation
of this argument.

\begin{proof}[Proof of Theorem \ref{thm:main:1}]
Note that it is enough to prove that the two sides of
(\ref{eq:mod:0}) agree in any local coordinate system. Hence,
let us choose local coordinates $(x^1,\dots,x^m)$, so that:
\begin{align*}
    \pi_0&=\sum_{i<j}\pi_0^{ij}\frac{\partial}{\partial
  x^i}\wedge\frac{\partial}{\partial x^j},\\
  \mc{N}&=\sum_{i,j}N^i_j\frac{\partial}{\partial x^i}\wedge\d x^j,\\
  \mu&=f\d x^1\wedge\cdots\wedge\d x^m.
\end{align*}
In these local coordinates, the compatibility condition
(\ref{eq:PN:compatibiliy:2}) for a Poisson-Nijenhuis structure reads:
\begin{multline*}
0=[\d x^i,\d x^j]_{\pi_1}-[\d x^i,\d x^j]_{\pi_0}^{\mc{N}^*}=\\
=\sum_{k,l}\left(\pi_0^{lj} \frac{\partial N_k^i}{\partial x^l}
+\pi_0^{il} \frac{\partial N_k^j}{\partial x^l}
-\pi_0^{lj} \frac{\partial N_l^i}{\partial x^k}
-N_k^l\frac{\partial \pi_0^{ij}}{\partial x^l}
+N_l^j\frac{\partial \pi_0^{il}}{\partial x^k}
\right)\d x^k
\end{multline*}
If in each coefficient of $\d x^k$ we contract $j$ and $k$, we see
that the two last terms cancel out, and we obtain:
\[ \sum_{k,l}\left(2\pi_0^{lk} \frac{\partial N_k^i}{\partial x^l}
+\pi_0^{il} \frac{\partial N_k^k}{\partial x^l}\right)=0, \qquad (i=1,\dots,m).\]
Using this identity and Proposition \ref{prop:mod:Poisson}, we conclude that:
\begin{align*}
X_{\mc{N}}
&=X^1_\mu-\mc{N} X^0_\mu\\
&=
\sum_{i,j} \frac{\partial \pi^{ij}_1}{\partial x_j} \frac{\partial}{\partial x_i}
-\mc{N}\sum_{i,j} \frac{\partial \pi^{ij}_0}{\partial x_j} \frac{\partial}{\partial x_i}\\
&=\sum_{i,j,k} \left(\frac{\partial \pi^{ik}_0 N^j_k}{\partial x_j}
-N^i_k\frac{\partial \pi^{kj}_0}{\partial x_j}\right)\frac{\partial}{\partial x_i}\\
&=\sum_{i,l,k}\pi_0^{lk} \frac{\partial N_k^i}{\partial x^l}\frac{\partial}{\partial x_i}\\
&=-\frac{1}{2}\sum_{i,l,k}\pi_0^{il} \frac{\partial N_k^k}{\partial x^l}\frac{\partial}{\partial x_i}
=X^0_{-\frac{1}{2}\tr \mc{N}}.
\end{align*}
\end{proof}

\begin{rem}
We recall (see \cite[page 58]{KosmannMagri}) that one has a
commutative diagram of morphisms of Lie algebroids:
\[
\xymatrix{
(T^*M,[\cdot,\cdot]_{\pi_1})\ar[rr]^{\mc{N}^*}\ar[dd]_{\pi^\sharp_0}\ar[ddrr]^{\pi^\sharp_1}
&&
(T^*M,[\cdot,\cdot]_{\pi_0})\ar[dd]^{\pi^\sharp_0}\\
\\
(TM,[\cdot,\cdot]_\mc{N})\ar[rr]^{\mc{N}}&&(TM,[\cdot,\cdot])
}
\]
The relative modular class of the morphism $\mc{N}$ on the bottom
horizontal arrow is represented by $\d(\tr\mc{N})$. On the other
hand, Theorem \ref{thm:main:1} states that the relative modular
class of the morphism $\mc{N}^*$ on the top arrow is represented by
$-(\pi_0^\sharp)\d(\tr\mc{N})=(\pi_0^\sharp)^*\d(\tr\mc{N})$.
Hence this diagram codifies nicely the relationship between all
(relative and absolute) modular classes involved in a Poisson-Nijenhuis
manifold.
\end{rem}

\section{Integrable hierarchies and the modular class}

We  now consider integrable hierarchies of hamiltonian systems and
observe that  they are closely related to the modular vector field
introduced above.

\subsection{The hierarchy of a non-degenerate PN manifold}
Theorem \ref{thm:main:1} above shows that for any non-degenerate
Poisson-Nijenhuis manifold the modular vector field $X_{\mc{N}}$
is hamiltonian relative to $\pi_0$ and Poisson relative to
$\pi_1$. So, the question arises whether this vector field is also
Hamiltonian  relative to $\pi_1$ with respect to another function
of $\NN$, i.e., whether it is a bi-hamiltonian vector field. Of
course, the obstruction is the Poisson cohomology class
$[X_{\mc{N}}]$, the modular class of the Poisson-Nijenhuis
manifold.  An important case where this class vanishes and
which fits many examples is the following:

\begin{theorem}
\label{thm:main:2}
Let $(M,\pi_0,\mc{N})$ be a Poisson-Nijenhuis manifold and assume that
$\mc{N}$ is non-degenerate. Then the modular vector field
$X_{\mc{N}}$ is bi-hamiltonian and hence determines a hierarchy of
flows which is given by:
\begin{equation}
\label{eq:hier:PN} X_{i+j}=\pi_i^\sharp\d h_j=\pi_j^\sharp\d
h_i,\quad (i,j \in {\bf Z})
\end{equation}
where
\begin{equation}
\label{eq:hamilt:PN} h_0=-\frac{1}{2}\,\log(\det\mc{N}),\quad
h_i=-\frac{1}{2i}\,\tr\mc{N}^i, \quad (i \not= 0).
\end{equation}
\end{theorem}

\begin{proof}
Let us start by verifying that $X_{\mc{N}}$ is bi-hamiltonian:
\begin{equation}
\label{eq:biham:PN}
X_{\mc{N}}=X^1_{-\frac{1}{2}\log(|\det\mc{N}|)}=X^0_{-\frac{1}{2}\tr
  \mc{N}},
\end{equation}
By Theorem \ref{thm:main:1}, we just need to prove the first equality.

We claim that the first equality holds on the open dense
set of common regular points of $\pi_0$ and $\pi_1$. In fact, for any such
regular point we can choose an open neighborhood $U$ where both $\pi_0$ and $\pi_1$
admit invariant volume forms $\mu_0$ and $\mu_1$. It is easy to see that we can
take these two volume forms to be related by $\mc{N}$:
\[ \mu_1=\mc{N}^{-1}\mu_0=\frac{1}{|\det(\mc{N})|^{\frac{1}{2}}}\mu_0.\]
where $n=\frac{1}{2}\dim M$. It follows from relation (\ref{eq:change:form})
that the modular vector fields for $\pi_1$ relative to these two n-forms are related by:
\[ 0=X^1_{\mu_1}=X^1_{\mu_0}+X^1_{\frac{1}{2}\log(|\det\mc{N}|)}.\]
Since $X^0_{\mu_0}=0$, we can compute the modular vector field $X_{\mc{N}}$
as follows:
\[
X_{\mc{N}}=X^1_{\mu_0}-\mc{N}X^0_{\mu_0}=X^1_{-\frac{1}{2}\log(|\det\mc{N}|)}.
\]
This proves our claim, so $X_{\mc{N}}$ is
bi-hamiltonian(\footnote{F.~Magri as pointed out to us that
relation (\ref{eq:biham:PN}) also follows from the identity
$\mc{N}^*\d(\det\mc{N}))=\det\mc{N}\d(\tr\mc{N})$, which a
consequence of the fact that $\mc{N}$ is a root of its
characteristic polynomial.})

Now, it remains to prove the multi-hamiltonian structure for the
higher flows. This follows by an iterative procedure. For example, let
us check the multi-hamiltonian structure of the 2nd flow in the hierarchy:
\begin{equation}
\label{eq:hier:2}
X_2=X^2_{-\frac{1}{2}\,\log(\det\mc{N})}
=X^1_{-\frac{1}{2}\,\tr\mc{N}}=X^0_{-\frac{1}{4}\,\tr\mc{N}^2}
\end{equation}
First, we note that, from what we just proved, we have:
\[
X_{\mc{N}^2}=X^2_{\mu}-\mc{N}^2X^0_{\mu}=X^0_{-\frac{1}{2}\,\tr\mc{N}^2}
=X^2_{-\frac{1}{2}\det(\log\NN^2)}.
\]
which shows equality of two terms in (\ref{eq:hier:2}). On the
other hand, we have:
\begin{align*}
X^1_{-\frac{1}{2}\,\tr\mc{N}}&=X^2_{\mu}-\mc{N}X^1_{\mu}\\
&=X^2_{\mu}-\mc{N}^2X^0_{\mu}+\mc{N}^2X^0_{\mu}-\mc{N}X^1_{\mu}\\
&=X_{\mc{N}^2}-NX_{\mc{N}}\\
&=X^0_{-\frac{1}{2}\,\tr\mc{N}^2}-NX^0_{-\frac{1}{2}\,\tr\mc{N}}
=X^0_{-\frac{1}{2}\,\tr\mc{N}^2}-X^1_{-\frac{1}{2}\,\tr\mc{N}}.
\end{align*}
which gives:
\[ X^1_{-\frac{1}{2}\,\tr\mc{N}}=X^0_{-\frac{1}{4}\,\tr\mc{N}^2},\]
so giving equality with the remaining term in (\ref{eq:hier:2}).

By iteration, looking at the vector fields $X_{\mc{N}^i}$, one
obtains the multi-hamiltonian formulation of the remaining higher
order flows. For negative values of the index we apply the
proposision to $\NN^{-1}$ to obtain

\[ X^1_{-\frac{1}{2}\log(|\det\mc{N}|)}=X^0_{\frac{1}{2}\tr
  \NN^{-1}}, \]
and then we proceed as in the case of positive indices.
\end{proof}

Note that $\mc{N}$ always has double eigenvalues. For the hierarchy to
be completely integrable, we need $n=\frac{1}{2}\dim M$ independent
spectral invariants $\det\NN, \tr\NN, \tr\NN^2, \dots$ This will
follow if $\mc{N}$ has $n=\frac{1}{2}\dim M$ independent eigenvalues.

\subsection{Master symmetries and modular vector fields}
\label{sec:Oevel}

When $\mc{N}$ is degenerate the results in the previous paragraph
do not apply. In this situation, there is a procedure due to Oevel
\cite{oevel2} to produce integrable hierarchies from master
symmetries, and it is natural to look how the modular vector
fields fit into this scheme. We start with the following result
which is of independent interest:

\begin{prop}
\label{prop:modular:def}
Let $\pi_0$ and $\pi_1$ be Poisson tensors such that
$\pi_1=\Lie_Z\pi_0$, for some vector field $Z$. Also fix a volume
form $\mu\in\Omega^{\text{top}}(M)$. Then their modular vector fields
are related by
\begin{equation}
\label{eq:modular:def}
X^1_\mu=\Lie_Z X^0_\mu+X^0_{\div_\mu(Z)}.
\end{equation}
\end{prop}

\begin{proof}
The proof is straightforward if we use the definition of the modular
class in terms of the homological operator $D_\mu$:
\begin{align*}
X^1_\mu&=D_\mu(\pi_1)\\
&=D_\mu([Z,\pi_0])\\
&=[Z,D_\mu(\pi_0)]-[D_\mu(Z),\pi_0]\\
&=\Lie_Z X^0_\mu+X^0_{\div_\mu(Z)}.
\end{align*}
\end{proof}

Now, in Oevel's approach, one assumes that we have a
bi-hamiltonian system defined by the Poisson tensors $\pi_0$ and
$\pi_1$ and the hamiltonians $h_1$ and $h_0$:
\begin{equation}
\label{eq:ini:biham}
X_1=X^0_{h_1}\equiv\pi^\sharp_0\d h_1=X^1_{h_0}\equiv\pi^\sharp_1\d
h_0.
\end{equation}
If, additionaly, $\pi_0$ is symplectic, one can define the recursion
operator in the usual way:
\[ \mc{N} = \pi^\sharp_1 \circ (\pi^\sharp_0)^{-1},\]
the higher flows $X_{i}:=\mc{N}^{i-1} X_1$, and the higher order
Poisson tensors $\pi_i:=\mc{N}^i \pi_0$. Note that $\mc{N}$ can now be
degenerate. Now one can generate master-symmetries by the following method:

\begin{theorem}[\cite{oevel2}]
Suppose that $Z_0$ is a conformal symmetry for both $\pi_0$, $\pi_1$
and $h_0$, i.e.,  for some scalars $\lambda$, $\mu$, and $\nu$ we
have
\[
\Lie_{Z_0}\pi_0= \lambda \pi_0,\quad
\Lie_{Z_0}\pi_1= \mu \pi_1, \quad
\Lie_{Z_0} h_0 = \nu h_0.
\]
Then the vector fields
\[ Z_i =\mc{N}^i Z_0 \]
are master symmetries  and  we have
\begin{align}
\Lie_{Z_i} h_j &= (\nu+(j-1+i)(\mu-\lambda)) h_{i+j},\\
\Lie_{Z_i}\pi_j&= (\mu+(j-i-2)(\mu-\lambda)) \pi_{i+j},\\
      [Z_i,Z_j]&= (\mu-\lambda)(j-i) Z_{i+j}.
\end{align}
\end{theorem}

To simplify the notation, we will set:
\[ c_{i,j}= (\mu +(j-i-2) (\mu -\lambda)).\]
so, for example, $[Z_i,\pi_j]=c_{i,j}\pi_{i+j}$. Also, we fix a volume
form $\mu\in\Omega^{\text{top}}(M)$, so the $j$th Poisson bracket
in the hierarchy has the modular vector field
\[ X^j_\mu=D_\mu(\pi_j).\]
The following proposition establishes relations among these modular vector
fields:

\begin{theorem}
\label{prop:mod:vf:hierarchy}
For the hierarchy above:
\begin{align*}
[X^j_\mu,Z_i]&=c_{i,j}X^{i+j}_\mu+X_{f_i}^j,\\
\Lie_{X^i_\mu}\pi_j&=-\Lie_{X^j_\mu}\pi_i,
\end{align*}
where $f_i=D_\mu(Z_i)=\div_\mu(Z_i)$.
\end{theorem}

\begin{proof}
To prove the first relation, one simply applies  Proposition
\ref{prop:modular:def} repeatedly. For the second relation, we
observe that:
\begin{align*}
\Lie_{X^i_\mu}\pi_j &=[X^i_\mu,\pi_j] \\
                &=\frac{1}{c_{j-i,i}}\left[X^i_\mu,[Z_{j-i},\pi_i]\right]
                \ .
\end{align*}
Using the super--Jacobi identity for the Schouten bracket and the
fact that $X^i_\mu$ is Poisson relative to $\pi_i$, the last term reduces to
$\left[Z_{j-i},X^i_\mu],\pi_i \right]$. Therefore
\begin{align*}
\Lie_{X^i_\mu}\pi_j
&=\frac{1}{c_{j-i,i}}\left[[Z_{j-i},X^i_\mu],\pi_i\right] \\
&=-\frac{1}{c_{j-i,i}}\left[[Z_i,X^{i-j}_\mu],\pi_i\right] \\
&=-\frac{1}{c_{j-i,i}}\left[c_{j-i,i}X^j_\mu+X_{f_{j-i}}^i,\pi_i\right] \\
&=-[X^j_\mu,\pi_i]=-\Lie_{X^j_\mu}\pi_i.
\end{align*}
\end{proof}

Note that, even when $\mc{N}$ is non-degenerate, there is no
reason for the hierarchy of flows produced by this method to
coincide with the hierarchy of flows canonically  associated with
the Poisson-Nijenhuis manifold. In general, one would obtain two
distinct hierarchies. However, as we shall see in the next
section, in most of the  examples it is often the case that this
two hierarchies coincide. This is due to the fact that, in many
examples, the initial bi-hamiltonian system (\ref{eq:ini:biham})
has a multiple of $\tr\mc{N}$ as one of the hamiltonians.

\section{Examples}
In this section we will illustrate the results of this paper on some
well-known integrable systems such as the Harmonic oscillator, the
Calogero-Moser system  and various versions of the Toda lattice.
This gives a new approach to the bi-hamiltonian structure of these
systems, and in some cases leads to some new results.

\subsection{Harmonic oscillator}
This classical integrable system has a well-known bi-hamiltonian
structure which we recall using the notation from
(\cite{caseiro}).

On $\Rr^{2n}$ with the standard symplectic structure and canonical
coordinates $(q_i, p_i)$, consider the following hamiltonian function:
\[ h_1=\sum_{i=1}^n \frac12 (p_i^2+q_i^2) \ . \]
The resulting hamiltonian system is completely integrable with the
following integrals of motion in involution:
\[I_i=\frac12 (p_i^2+q_i^2), \quad (i=1,\dots,n). \]
For its bi-hamiltonian structure one takes the Poisson structure
associated with the canonical symplectic form:
\[ \pi_0=\sum_{i=1}^n \frac {\partial }{\partial p_i} \wedge  \frac
{\partial} {\partial q_i}, \]
and the new Poisson structure:
\[ \pi_1=\sum_{i=1}^n  I_i \frac {\partial }{\partial p_i} \wedge  \frac
{\partial} {\partial q_i} \ . \]
These form a compatible pair of Poisson structures, and we also have:
\begin{equation}
\label{eq:hamornic:osc}
X_1=\sum_{i=1}^n  p_i \frac {\partial }{\partial q_i} - q_i \frac
{\partial} {\partial p_i}
=\pi_0^\sharp\d h_1=\pi_1^\sharp \d h_0,
\end{equation}
where
\[ h_0=\log{I_1}+ \dots + \log{I_n}.\]

It is easy to see that this is the bi-hamiltonian formulation of
the first flow in the integrable hierarchy of the
Poisson-Nijenhuis manifold $(\pi_0,\mc{N})$, where the
Nijenhuis tensor is the diagonal (1,1) tensor:
\[ \mc{N}=\text{diag}(I_1,\dots,I_n,I_1,\dots,I_n). \]
In fact, with this definition, we find $\pi_1=\NN\pi_0$ and:
\[ \det \NN = \prod_{i=1}^n I_i^2, \]
so that:
\begin{align*}
{\frac12}\log(\det\NN)&=\log I_1 + \cdots + \log I_n=h_0,\\
{\frac12}\tr\NN&=I_1+\cdots+I_n=h_1.
\end{align*}
Hence, the bi-hamiltonian formulation (\ref{eq:hamornic:osc})
coincides with the one of the first flow of the hierarchy
(\ref{eq:hier:PN}) canonically  associated with the
Poisson-Nijenhuis manifold $(\pi_0,\mc{N})$.

In this example, we have a mastersymmetry $Z$ such that
$\Lie_Z\pi_0=\pi_1$ which is given by:
\[ Z=-\sum_{i=1}^n \frac14 I_i \left(q_i\frac{\partial}{\partial q_i}
+ p_i\frac{\partial}{\partial p_i}\right).\]
If we let $\mu_0=\d p_1\wedge\d q_1\wedge\cdots\wedge\d p_n\wedge\d q_n$,
which is a Liouville form for $\pi_0$, we compute:
\[ -\div_{\mu_0}(Z)=\frac12\, \tr \NN=h_1,\]
as expected.

There is however one point that we overlooked: strictly speaking
these results are true only on the manifold
$M=\Rr^{2n}-\cup_{i=1}^n\{I_i=0\}$ where $\mc{N}$ is invertible.
In fact, on $\Rr^{2n}$ the relative modular vector field $X_\NN$
is not hamiltonian relative to $\pi_1$! In case
$X_\NN=\pi_1^\sharp\d H$ for some smooth function $H\in
C^\infty(\Rr^{2n})$ then, on points away from $I_i=0$, $H$ must
differ from $h_0=\log{I_1}+ \dots + \log{I_n}$ by a constant, and
this is clearly impossible. Therefore, on $\Rr^{2n}$ the relative
modular class $[X_\NN]$ is non-trivial, and there is no canonical
bi-hamiltonian hierarchy.

Note that this examples is \emph{universal}: any integrable hierarchy
associated with a non-degenerate Poisson-Nijenhuis manifold
$(M,\pi_0,\mc{N})$ locally (in action-variables coordinates) looks
like this one.

\subsection{The rational Calogero-Moser system}

The Calogero-Moser system is a well-known finite-dimensional
integrable system (in fact, super-integrable). One can define this
system on $\Rr^{2n}$, with the standard symplectic structure and
canonical coordinates $(q_i,p_i)$, by the hamiltonian function:
\[
h_2=\frac12 \sum_{i=1}^n p_i^2 +
\frac{g^2}{2}\sum_{j\not=i}\frac{1}{(q_i-q_j)^2}.
\]

The Calogero-Moser system admits a Lax pair formulation where the Lax
matrix $L$ is given by
\[ L_{ij}=p_i\delta_{ij}+g\frac{i(1-\delta_{ij})}{q_i-q_j}. \]
The system is then completely integrable with involutive first integrals
given by:
\[ F_i=\tr (L^i),\quad (i=1,\dots,n). \]
 Moreover, following Ranada \cite{ranada}, consider also the
functions $G_i=\tr(Q L^{i-1})$, where $Q$ is the diagonal matrix
$\text{diag}(q_1,\dots,q_n)$. It was shown in \cite{francoise}
that these functions are independent and lead to the algebraic
linearization of the system. Using these functions as coordinates,
we can write the hamiltonian vector field in the form:
\[ X_1=\sum_{i=1}^n F_i\frac{\partial}{\partial G_i}. \]
The original Poisson structure becomes:
\[ \pi_0=\sum_{i=1}^n \frac {\partial }{\partial F_i} \wedge  \frac
{\partial} {\partial G_i},\]
and there exists a second compatible Poisson structure given by:
\[ \pi_1=\sum_{i=1}^n  F_i \frac {\partial }{\partial F_i}\wedge
\frac{\partial}{\partial G_i},\]
providing a bi-hamiltonian formulation given by
\[
X_1=\pi_0^\sharp\d h_i=\pi_1^\sharp \d h_{i-1}\quad (i=2,\dots,n).
\]
where
\[
h_j=\frac {1}{2j} \tr\left(\pi_1^\sharp\circ(\pi_0^\sharp)^{-1}\right)^j
=\frac{1}{2j}\sum_k \left(F_k \right)^j, \quad  (j=1,\dots,n).
\]
Now we observe that if we let $h_0:=\log(F_1\cdots F_n)$,
then we can write the system in the form:
\[ X_1=\pi_0^\sharp \d h_1 = \pi_1^\sharp \d h_0.\]
If we set $\NN:=\pi_1^\sharp\circ(\pi_0^\sharp)^{-1}$, then one checks
easily that:
\[ h_0=\frac{1}{2}\,\log(\det \NN),\quad h_1=\frac{1}{2}\tr \NN, \]
so $X_1$ is in fact the first flow of the hierarchy
(\ref{eq:hier:PN}) canonically associated with the
Poisson-Nijenhuis manifold $(\Rr^{2n},\pi_0,\mc{N})$.

\subsection{Toda lattice in Moser coordinates}
Our next example is related to the Toda hierarchy in the so-called
Moser coordinates. The hierarchy of Poisson tensors is due to
Faybusovich and Gekhtman \cite{fay}, and can be defined as
follows. Consider $\Rr^{2n}$ with coordinates
$(\lambda_1,\dots,\lambda_n,r_1,\dots,r_n)$ and define the Poisson
structures:
\begin{align*}
\pi_0&=\sum_{i=1}^n r_i
\frac{\partial}{\partial\lambda_i}\wedge \frac{\partial}{\partial r_i},\\
\pi_1&=\sum_{i=1}^n \lambda_i r_i
\frac{\partial}{\partial\lambda_i}\wedge \frac{\partial}{\partial r_i}.
\end{align*}
We also set $h_2=\frac{1}{2} \sum_{j=0}^n \la_j^2$ as the
hamiltonian. Using the $\pi_0$ bracket, we obtain the following set of
Hamilton's equations
\[ \dot{\la}_i=0, \qquad \dot{r}_i=\la_i r_i, \qquad (i=1,\dots,n). \]
This system is bi-hamiltonian, since
\[\pi_0^\sharp \d h_2= \pi_1^\sharp \d h_1 \ , \]
where $h_1=\sum_{j=1}^n \la_j$.
Again, if we define
\[ \NN:=\pi_1^\sharp\circ(\pi_0^\sharp)^{-1}=
\text{diag}\left(\la_1,\dots,\la_n,\la_1,\dots,\la_n \right), \]
we have:
\[
\tr \NN=2\, \sum_{j=1}^n \la_j=2\, h_1,\qquad
\tr \NN^2=4\, \sum_{j=1}^n \la_j^2=4\, h_2.
\]
It follows that our system is in fact the \emph{second} flow in
the hierarchy (\ref{eq:hier:PN}) canonically  associated with the
Poisson-Nijenhuis manifold $(\Rr^{2n},\pi_0,\mc{N})$. The first
flow of this hierarchy is (\footnote{Note
that, just like in the case of the harmonic oscilator we should
exclude the points with some $\la_i=0$, where $\det\NN$ vanishes.}):
\[ X_1=\pi_0^\sharp \d h_1 = \pi_1^\sharp \d h_0.\]
where $h_0=\frac{1}{2}\log(\det\NN)=\log\la_1+\cdots+\log\la_n$, and
in coordinates is simply:
\[ \dot{\la}_i=0, \qquad \dot{r}_i=r_i, \qquad (i=1,\dots,n). \]

In this example, there is a mastersymmetry connecting $\pi_0$ and
$\pi_1$ given by:
\[ Z=-\frac{1}{2}\sum_{j=1}^n \la_j^2 \frac{\p}{\p\la_j}\]
so that $\Lie_{Z}\pi_0=\pi_1$. Then, as expected, we find:
 \[ \div(Z)=-\sum_{j=1}^n \la_j=-h_1. \]

This example also falls in Oevel's scheme of Section \ref{sec:Oevel}.
The vector field:
\[ Z_0=\sum_{j=1}^n \la_j \frac{\p}{\p\la_j},\]
is a conformal symmetry of $\pi_0$, $\pi_1$ and $h_1$:
\[ \Lie_{Z_0}\pi_0=-\pi_0,\quad \Lie_{Z_0}\pi_1=0, \quad \Lie_{Z_0}h_1=h_1.\]
Note here $h_1$, instead of $h_0$. This means that it is the second flow in
the hierarchy (i.e., the original flow) that falls into Oevel scheme!
Recalling now that $Z_i=\NN^i Z_0$, we find:
\[ Z_{-1}=\sum_{j=1}^n \frac{\p} { \p \la_j},\qquad Z_1=-2\, Z.\]
If we let $\mu$ be the standard volume on $\Rr^{2n}$, we see that
$Z_{-1}$ coincides with the modular vector field for $\pi_0$ relative
to $\mu$:
\[ X_{\mu}^0=\sum_{j=1}^n \frac{\p}{\p \la_j}=Z_{-1}. \]
On the other hand
\[ X_{\mu}^1=\sum_{j=1}^n \la_j\frac{\p}{\p \la_j}-
\sum_{j=1}^n r_j\frac {\p}{ \p r_j},\]
so that:
\[
\X_{\mc{N}}=X_{\mu}^1-\mc{N} X_{\mu}^0
=X_{\mu}^1-\mc{N} Z_{-1}
=X_{\mu}^1-Z_0=X^1_{h_0},
\]
is indeed the first flow in the hierarchy.

\subsection{Bogoyavlensky-Toda systems}
We consider now the example of the $C_n$ Toda system. The $B_n$ and
$D_n$ Toda systems are similar and details on the computations
can be found in \cite{damianou3,damianou4}.

{To} define the $C_n$ system one considers $\Rr^{2n}$ with the
canonical symplectic structure and the hamiltonian function:
\[
H_2=\frac{1}{2}\,\sum_1^n p_j^2 + e^{ q_1- q_2} + \cdots + e^{ q_{n-1}-q_n}
   + e^{ 2 q_n}.
\]
Let us consider the Flaschka-type change of coordinates:
\begin{align*}
  a_i&=\frac{1}{2}\,e^{\frac{1}{2}(q_i-q_{i+1})},\qquad (i=1,\dots,n-1)\\
  a_n&=\frac{1}{\sqrt{2}}\,e^{q_n},\\
  b_i&=-\frac{1}{2}\,p_i,\qquad (i=1,\dots,n).
\end{align*}
The equations for the flow in $(a_i,b_i)$ coordinates become
\begin{align*}
 \dot a_i&=a_i\,(b_{i+1}-b_i ),\qquad (i=1,\dots,n-1),\\
 \dot a_n&=-2\,a_n b_n,\\
 \dot b_i&= 2\,(a_i^2-a_{i-1}^2 ),\qquad (i=1,\dots,n),
\end{align*}
with the convention that $a_0=0$. These equations can also be written
as a Lax pair  $\dot L=[B,L]$, where the Lax matrix $L$ is given by:
\[
L=  \begin{pmatrix} b_1 &  a_1 &  & &  &    &    & \cr
                   a_1 &  \ddots  & \ddots  &   & &&   &  \cr
                    & \ddots & \ddots& a_{n-1} &  & &   &  \cr
                   &  & a_{n-1} & b_n & a_n &  & &  \cr
                    & &  & a_n & -b_n &  -a_{n-1}     & & \cr
                    & & & & -a_{n-1} & \ddots & \ddots &  \cr
                    &  &&&& \ddots & \ddots & -a_1  \cr
                       &&&&& & -a_1 & -b_1     \end{pmatrix}   \ ,
\]
and $B$ is the skew-symmetric part of $L$.

In the new variables $(a_i,b_i)$, the canonical Poisson bracket on
$\Rr^{2n}$ is transformed into a bracket $\pi_1$ which is given by
\begin{align*}
\{ a_i, b_i \}&=-a_i,\qquad (i=1,\dots,n-1)\\
\{ a_i, b_{i+1}\}&=a_i,\qquad (i=1,\dots,n-1)\\
\{a_n, b_n\}&=-2 a_n.
\end{align*}
We follow the tradition of denoting this bracket by $\pi_1$ (instead
of $\pi_0$) being a linear bracket of degree one. This will lead to a
shift in degrees, when compared to the formulas in the rest of the
paper (to obtain the same formulas we should denote this bracket by
$\pi_0$). The same comments applies to the first integrals of the
system which, following the tradition, will be denoted by $H_2,
H_4,\dots, H_{2n}$, where:
\[ H_{2i} = \frac{1}{2i} \tr L^{2i}.\]

In order to obtain a bi-hamiltonian formulation (see \cite{joana}),
one introduces a cubic Poisson bracket $\pi_3$, defined by:
\begin{align*}
\{a_i, a_{i+1}\}&=a_i a_{i+1} b_{i+1},\qquad (i=1,\dots,n-2)\\
\{a_{n-1},a_n \}&=2\,a_{n-1}a_n b_n, \\
\{ a_i, b_i \}  &=-a_i b_i^2 -a_i^3,\qquad (i=1,\dots,n-1)\\
\{ a_n,b_n \}   &=-2\,a_n b_n^2 -2\,a_n^3,\\
\{a_i, b_{i+2}\}&=a_i a_{i+1}^2,\qquad\qquad (i=1,\dots,n-2)\\
\{ a_i,b_{i+1}\}&=a_i b_{i+1}^2 + a_i^3,\qquad (i=1,\dots,n-1)\\
\{ a_{n-1},b_n\}&=a_{n-1}^3+a_{n-1}(b_n^2-a_n^2),\\
\{ a_i,b_{i-1}\}&=-a_{i-1}^2 a_i,  \qquad\qquad (i=1,\dots,n)\\
\{ a_n,b_{n-1}\}&=-2\,a_{n-1}^2  a_n, \\
\{ b_i,b_{i+1}\}&=2\,a_i^2(b_i+b_{i+1}),\qquad (i=1,\dots,n-1).
\end{align*}
This leads immediately (\cite{joana}) to a bi-hamiltonian system:
\[ \pi_3^\sharp \d h_2 = \pi_1^\sharp \d h_4.\]
However, this is not the original system. For the original system, we
follow \cite{damianou3} and define $\pi_{-1}=\pi_1 \pi_3^{-1} \pi_1$.
Then the $C_n$ Toda system has the bi-hamiltonian formulation:
\[ \pi_1^\sharp \d h_2 = \pi_{-1}^\sharp \d h_4.\]
We now give a new bi-hamiltonian formulation using our Theorem
\ref{thm:main:2}. We have the Nijenhuis tensor:
\[ \NN:=\pi_3^\sharp\circ(\pi_1^\sharp)^{-1},\]
and we set
\[ H_0:=\frac{1}{2}\log(\det\NN). \]
We need to check that the hamiltonian vector field of
$H_0$ with respect to the second bracket $\pi_3$ satisfies:
\begin{equation}
\label{lenard}
\pi_3^\sharp \d H_0=\pi_1^\sharp \d H_2,
\end{equation}
so that this yields a bi-hamiltonian formulation for the $C_n$-Toda.
In fact, this follows easily from the Lenard relations for the
eigenvalues of the Lax matrix $L$:
\[
\pi_3^\sharp \d \lambda_i= \lambda_i^2 \pi_1^\sharp \d  \lambda_i,
\]
which lead to:
\begin{align*}
\pi_3^\sharp \d H_0
          &=\pi_3^\sharp \d\frac 12 \log(\det L) \\
          &=\frac12 \sum_{1=1}^n \frac{1}{\lambda_I}\pi_3^\sharp\d\lambda_i \\
          &=\frac12 \sum_{i=1}^n \lambda_i \pi_1^\sharp \d\lambda_i \\
          &=\sum_{i=1}^n \lambda_i^2 \pi_1^\sharp \d \frac{\lambda_i^2}{2}
          =\pi_1^\sharp \d H_2.
\end{align*}
Of course, this computation can be avoided by invoking Theorem
\ref{thm:main:2}.

The situation in the other simple Lie algebras of type $B_N$ and
$D_N$ is entirely similar. Therefore we have the following result:

\begin{theorem}
Consider the  $B_n$, $C_n$ and $D_n$ Toda systems. In each case we
define
\[ \NN:=\pi_3^\sharp\circ(\pi_1^\sharp)^{-1},\]
where $\pi_1$ is the Lie-Poisson bracket and $\pi_{3}$ is the
cubic Poisson bracket. Also, let $H_0=\log(\det(L))$ and
$H_{2i}=\frac{1}{2i}\tr L^{2i}$, if $i \not=0$. Then we have the
following new bi-hamiltonian formulation for these systems:
\[
\pi_3^\sharp \d H_0=\pi_1^\sharp \d H_2
\quad \Longleftrightarrow\quad
X_{\tr\NN}^1=X_{\log{\det\NN}}^3.\]
The function $\sqrt{\det \NN}$ equals the determinant of
$L$ for $C_N$ and $D_N$ and the product of the non-zero eigenvalues of
$L$ for $B_N$, while the function $\frac12 \tr\NN=H_2$ is the
original hamiltonian. Finally,
\begin{equation}
\label{todachain}
\pi_{2k+1}^\sharp \d H_{2-2k}=\pi_{2k-1}^\sharp \d H_{4-2k},
\qquad (k\in\Zz).
\end{equation}
\end{theorem}

\subsection{Finite, non-periodic Toda lattice}
The case of the $A_n$ Toda lattice was already considered
in \cite{agrotis}, using specific properties of this system. We
use our general approach to show how one can quickly recover those
results.

The hamiltonian defining the Toda lattice is given in canonical
coordinates $(p_i,q_i)$ of $\Rr^{2n}$ by
\begin{equation}
\label{eq:Toda:An}
h_2(q_1,\dots,q_n,p_1,\dots,p_n)=
\sum_{i=1}^n\,\frac{1}{2}\,p_i^2+\sum _{i=1}^{n-1}\,e^{q_i-q_{i+1}}.
\end{equation}
For the integrability of the system we refer to the classical paper of
Flaschka \cite{flaschka1}.

Let us recall the bi-hamiltonian structure given in
\cite{fernandes1}. The first Poisson tensor in the hierarchy
is the standard canonical symplectic tensor, which we denote by
$\pi_0$, and the second Poisson tensor is:
\[
\pi_1=\begin{pmatrix} A_{n} & -B_{n} \\B_{n} & C_{n}
\end{pmatrix},
\]
where $A_n$, $B_n$ and $C_n$ are $n\times n$ skew-symmetric matrices
defined by
\begin{align*}
a_{ij}&=1=-a_{ji},\qquad (i<j)\\
b_{ij}&=p_i\delta_{ij},\\
c_{i,j}&=e^{q_i-q_{i+1}}\delta_{i,j+1}=-c_{j,i}, \qquad (i<j).
\end{align*}
Then setting $h_1=2\,(p_1+p_2+\cdots+p_n)$, we obtain the
bi-hamiltonian formulation:
\[ \pi_0^\sharp\d h_2=\pi_1^\sharp \d h_1.\]
If we set, as usual,
\[ \NN:=\pi_1^\sharp\circ(\pi_0^\sharp)^{-1},\]
then a small computation shows that Theorem \ref{thm:main:2} gives
the following multi-hamiltonian formulation:

\begin{prop}
The $A_n$ Toda hierarchy admits the multi-hamiltonian formulation:
\[ \pi_j^\sharp \d h_2=\pi_{j+2}^\sharp \d h_0,\]
where $h_0=\frac{1}{2}\log({\det \NN})$ and $h_2$ is the original
hamiltonian (\ref{eq:Toda:An}).
\end{prop}

If we change to Flaschka coordinates,
$(a_1,\dots,a_{n-1},b_1,\dots,b_n)$, then there is no recursion
operator anymore (recall that this is a singular change of
coordinates, where we loose one degree of freedom). Nevertheless, the
multi-hamiltonian structure does reduce (\cite{fernandes1}). One can
then compute the modular vector fields of the reduced Poisson tensors
$\pi_j$ relative to the standard volume form:
\[
\mu=\d a_1\wedge\cdots \wedge\d a_{n-1}\wedge\d b_1\wedge\cdots\wedge\d b_n.
\]
It turns out that the modular vector fields $X^j_\mu$ are hamiltonian
vector fields with hamiltonian function
\[ h=\log(a_1\cdots a_{n-1})+ (j-1) \log(\det(L)), \]
where $L$ is the Lax matrix.  For a discussion of this result we
refer to  \cite{agrotis}.  Note that the analogue of
(\ref{todachain}) in this case of the Toda chain is
\[ \pi_{j}^\sharp \d h_{2-j}=\pi_{j-1}^\sharp \d h_{3-j},   \quad k\in {\bf Z} \]
where $h_j=\frac{1}{j} \tr L^j$ for $j\not=0$ and
$h_0=\ln(\det(L))$.

\end{document}